\documentclass[12pt]{article}
\usepackage{amssymb}
\usepackage{amsmath}

\begin{document}

\vspace{.2in}\parindent=0mm

\begin{center}

 {\bf\Large {Pair of Dual Wavelet Frames  on Local Fields }}

\vspace{.5in}\parindent=0mm

{\bf{  M. Younus Bhat  }}

\parindent=0mm \vspace{.1in}

\end{center}

\parindent=0mm \vspace{.1in}
{{\it Department of  Mathematics,  National Institute of Technology, Srinagar-190001, Jammu and Kashmir, India. E-mail: $\text{gyounusg@gmail.com}$}}

\parindent=0mm \vspace{.2in}

{\bf{Abstract}} In this paper, an algorithm based on polyphase matrix for constructing a pair of orthogonal wavelet frames is suggested, and a general form for all orthogonal tight wavelet frames on local fields of positive characteristic is described. Moreover, we investigate their properties by means of the Fourier transform.
 \parindent=0mm \vspace{.2in}


\parindent=0mm \vspace{.2in}


\parindent=0mm \vspace{.2in}
{\bf{1. Introduction }}

\parindent=0mm \vspace{.2in}
An important example about frame is wavelet frame, which is obtained by translating and dilating a finite family of functions. A wavelet frame is a generalization of an orthonormal wavelet basis by introducing redundancy into a wavelet system. By sacrificing orthonormality and allowing redundancy, wavelet frames become much easier to construct than the orthonormal wavelets. Wavelet frames have many properties that make them useful in the study of function spaces, signal and image processing, sampling theory, optics, filter banks, wireless communications and so forth. In order to have more applications
of wavelet frames, several notions generalizing the concept of wavelet frames have been introduced and studied, namely dual wavelet frames [6], tight wavelet frames [7,8], pseudo wavelet frames [10] and orthogonal wavelet frames [2,3,9].

\parindent=0mm \vspace{.2in}
One of the most useful methods to construct tight wavelet frames is through the concept of unitary extension principle (UEP) introduced by Ron and Shen [11] and were subsequently extended by Daubechies et al.[4] in the form of the Oblique Extension Principle (OEP). They give sufficient conditions for constructing tight and dual wavelet frames for any given refinable function $\phi(x)$ which generates a multiresolution analysis. The resulting wavelet frames are based on multiresolution analysis, and the generators are often called framelets. Recent results in this direction can also be found in [1,5, 12-16] and the references therein.

\parindent=0mm \vspace{.2in}
Drawing inspiration from the construction of tight wavelet frames, we shall introduce the notion of orthogonal wavelet frames  on local fields of positive characteristic using extension principles. We present an algorithm for the construction of a pair of orthogonal wavelet frames based on polyphase matrices formed by the polyphase components of the wavelet masks. Moreover, we also gave a general construction algorithm for all orthogonal wavelet tight frames on local fields of positive characteristic from a compactly supported scaling function and investigate their properties by means of the Fourier transform.

\parindent=0mm \vspace{.2in}
The paper is structured as follows. In Section 2, we introduce some notations and preliminaries on local fields of positive characteristic including the definitions of Fourier transform and MRA based wavelet frame. In Section 3, we construct a pair of orthogonal wavelet frames and establish more conditions for the existence of orthogonal wavelet
frames in $L^2(R^+)$.

\parindent=0mm \vspace{.2in}
{\bf{1. Preliminaries on Local Fields }}

Let $K$ be a field and a topological space. Then $K$ is called a {\it local field} if both  $K^+$ and $K^*$ are locally compact Abelian groups, where  $K^+$ and $K^*$ denote the additive and multiplicative groups of $K$, respectively. If $K$ is any field and is endowed with the discrete topology, then $K$ is a local field. Further, if $K$ is connected, then $K$ is either $\mathbb R$ or $\mathbb C$. If $K$ is not connected, then it is totally disconnected. Hence by a local field, we mean a field $K$ which is locally compact, non-discrete and totally disconnected. The $p$-adic fields are examples of local fields. More details are referred to [13, 20]. In the rest of this paper, we use the symbols $\mathbb N, \mathbb N_0$ and $\mathbb Z$ to denote the sets of natural, non-negative integers and integers, respectively.

\parindent=8mm \vspace{.2in}
Let $K$ be a local field. Let $dx$ be the Haar measure on the locally compact Abelian group $K^{+}$. If $\alpha\in K$ and $\alpha\ne 0$, then $d(\alpha x)$ is also a Haar measure. Let $d(\alpha x)=|\alpha|dx$. We call $|\alpha|$ the {\it absolute value} of $\alpha$. Moreover, the  map $x\to |x|$ has the following properties: (a) $|x|=0$ if and only if $x = 0;$ (b) $|xy|=|x||y|$ for all $x, y \in K$; and (c) $|x+y|\le \max \left\{ |x|, |y|\right\}$ for all $x, y\in K$. Property (c) is called the {\it ultrametric inequality.} The set ${\mathfrak D}= \left\{x \in K: |x| \le 1\right\}$ is called the {\it ring of integers} in $K.$ Define ${\mathfrak B}= \left\{x \in K: |x| < 1\right\}$. The set ${\mathfrak B}$ is called the {\it prime ideal} in $K$. The prime ideal in $K$ is the unique maximal ideal in ${\mathfrak D}$ and hence as result ${\mathfrak B}$ is both principal and prime. Since the local field $K$ is totally disconnected, so there exist an element of ${\mathfrak B}$ of maximal absolute value. Let $\mathfrak p$ be a fixed element of maximum absolute value in ${\mathfrak B}$. Such an element is called a {\it prime element} of $K.$ Therefore, for such an ideal ${\mathfrak B}$ in ${\mathfrak D}$, we have ${\mathfrak B}= \langle \mathfrak p \rangle=\mathfrak p {\mathfrak D}.$ As it was proved in [20], the set ${\mathfrak D}$ is compact and open. Hence, ${\mathfrak B}$ is compact and open. Therefore, the residue space ${\mathfrak D}/{\mathfrak B}$ is isomorphic to a finite field $GF(q)$, where $q = p^{k}$ for some prime $p$ and $k\in\mathbb N$.

\parindent=8mm \vspace{.2in}
Let ${\mathfrak D}^*= {\mathfrak D}\setminus {\mathfrak B }=\left\{x\in K: |x|=1   \right\}$. Then, it can be proved that ${\mathfrak D}^*$ is a group of units in $K^*$ and if $x\not=0$, then we may write $x=\mathfrak p^k x^\prime, x^\prime\in {\mathfrak D}^*.$ For a proof of this fact we refer to [20]. Moreover, each ${\mathfrak B}^k= \mathfrak p^k {\mathfrak D}=\left\{x \in K: |x| < q^{-k}\right\}$ is a compact subgroup of $K^+$ and usually  known as the {\it fractional ideals} of $K^+$. Let ${\cal U}= \left\{a_i \right\}_{ i=0}^{q-1}$ be any fixed full set of coset representatives of ${\mathfrak B}$ in ${\mathfrak D}$, then every element $x\in K$ can be expressed uniquely  as $x=\sum_{\ell=k}^{\infty} c_\ell \mathfrak p^\ell $ with $c_\ell \in {\cal U}.$ Let $\chi$ be a fixed character on $K^+$ that is trivial on ${\mathfrak D}$ but is non-trivial on  ${\mathfrak B}^{-1}$. Therefore, $\chi$ is constant on cosets of ${\mathfrak D}$ so if $y \in {\mathfrak B}^k$, then $\chi_y(x)=\chi(yx), x\in K.$ Suppose that $\chi_u$ is any character on $K^+$, then clearly the restriction $\chi_u|{\mathfrak D}$ is also a character on ${\mathfrak D}$. Therefore, if $\left\{u(n): n\in\mathbb N_0\right\}$ is a complete list of distinct coset representative of ${\mathfrak D}$ in $K^+$, then, as it was proved in [13, 20], the set  $\left\{\chi_{u(n)}: n\in\mathbb N_0\right\}$   of distinct characters on ${\mathfrak D}$ is a complete orthonormal system on ${\mathfrak D}$.

\parindent=8mm \vspace{.2in}
The Fourier transform $\hat f$ of a function $f \in L^1(K)\cap L^2(K)$ is defined by\\
$$\hat f(\xi)= \displaystyle \int_K f(x)\overline{ \chi_\xi(x)}dx.\eqno(2.1)$$
It is noted that\\
$$\hat f(\xi)=  \int_K f(x)\,\overline{ \chi_\xi(x)}dx=  \int_K f(x)\chi(-\xi x)dx.$$

\parindent=0mm \vspace{.2in}
Furthermore, the properties of Fourier transform on local field $K$ are much similar to those of on the real line. In particular Fourier transform is unitary on $L^2(K)$.

\parindent=8mm \vspace{.2in}
We now impose a natural order on the sequence $\{u(n)\}_{n=0}^\infty$. We have ${\mathfrak D}/ \mathfrak B \cong GF(q) $ where $GF(q)$ is a $c$-dimensional vector space over the field $GF(p)$. We choose a set $\left\{1=\zeta_0,\zeta_1,\zeta_2,\dots,\zeta_{c-1}\right\}\subset {\mathfrak D^*}$ such that span  $\left\{\zeta_j\right\}_{j=0}^{c-1}\cong GF(q)$. For $n \in \mathbb N_0$ satisfying
$$0\leq n<q,~~n=a_0+a_1p+\dots+a_{c-1}p^{c-1},~~0\leq a_k<p,~~\text{and}~k=0,1,\dots,c-1,$$

\parindent=0mm \vspace{.1in}
we define
$$u(n)=\left(a_0+a_1\zeta_1+\dots+a_{c-1}\zeta_{c-1}\right){\mathfrak p}^{-1}.\eqno(2.2)$$

\parindent=0mm \vspace{.1in}
Also, for $n=b_0+b_1q+b_2q^2+\dots+b_sq^s, ~n\in \mathbb N_{0},~0\leq b_k<q,k=0,1,2,\dots,s$, we set

$$u(n)=u(b_0)+u(b_1){\mathfrak p}^{-1}+\dots+u(b_s){\mathfrak p}^{-s}.\eqno(2.3)$$

\parindent=0mm \vspace{.1in}
This defines $u(n)$ for all $n\in \mathbb N_{0}$. In general, it is not true that $u(m + n)=u(m)+u(n)$. But, if $r,k\in\mathbb N_{0}\; \text{and}\;0\le s<q^k$, then $u(rq^k+s)=u(r){\mathfrak p}^{-k}+u(s).$ Further, it is also easy to verify that $u(n)=0$ if and only if $n=0$ and $\{u(\ell)+u(k):k \in \mathbb N_0\}=\{u(k):k \in \mathbb N_0\}$ for a fixed $\ell \in \mathbb N_0.$ Hereafter we use the notation $\chi_n=\chi_{u(n)}, \, n\ge 0$.

\parindent=8mm \vspace{.2in}
Let the local field $K$ be of characteristic $t>0$ and $\zeta_0,\zeta_1,\zeta_2,\dots,\zeta_{c-1}$ be as above. We define a character $\chi$ on $K$ as follows:
$$\chi(\zeta_\mu {\mathfrak p}^{-j})= \left\{
\begin{array}{lcl}
\exp(2\pi i/t),&&\mu=0\;\text{and}\;j=1,\\
1,&&\mu=1,\dots,c-1\;\text{or}\;j \neq 1.
\end{array}
\right. \eqno(2.4)
$$
\parindent=0mm \vspace{.1in}
For given $\Psi := \left\{\psi_1,\dots, \psi_L\right\}\subset  L^2(K)$, define the  wavelet system

$${\mathcal F}(\Psi)=\Big\{ \psi_{\ell,j,k}(x)=q^{j/2}\psi_{\ell}\big({\mathfrak p}^{-j}x-u(k)\big),\, j\in \mathbb Z,k\in \mathbb N_0, \ell=1,2,\dots,  L\Big\}.\eqno(2.5)$$

\parindent=0mm \vspace{.1in}
The wavelet system ${\mathcal F}(\Psi)$ is called a  {\it  wavelet frame}, if there exist positive constants $A$ and $B$ such that

$$A\big\|f \big\|^2_{2} \le \sum_{\ell=1}^{L}\sum_{j\in\mathbb Z}\sum_{k\in \mathbb N_0} \left|\big\langle f, \psi_{\ell,j, k}\big\rangle\right|^2 \le B \big\|f\big\|^2_{2},\eqno(2.6)$$

\parindent=0mm \vspace{.1in}
holds for every $f\in  L^2(K)$, and we call the optimal constants $A$ and $B$ the lower frame bound and the upper frame bound, respectively. A {\it tight wavelet frame} refers to the case when $A = B$, and a Parseval wavelet frame refers to the case when $A = B = 1$. On the other hand if only the right hand side of the above double inequality holds, then we say ${\cal F}(\Psi)$ a {\it Bessel system}.

\parindent=8mm \vspace{.2in}
 Corresponding to the system (2.5), we have the dual system as

$${\cal F}( \Phi)=\Big\{\phi^{\ell}_{j,k}:=q^{j/2}\phi_{\ell}\big({\mathfrak p}^{-j}x-u(k)\big),\, j\in \mathbb Z,k\in \mathbb N_0, \ell=1,2,\dots,  L\Big\}.\eqno(2.7)$$

\parindent=0mm \vspace{.1in}
If both ${\mathcal F}(\Psi)$ and ${\cal F}( \Phi)$ are wavelet frames and for any $f\in L^2(K)$,  we have the reconstruction formula
$$f=\sum_{\ell=1}^{L}\sum_{j\in\mathbb Z}\sum_{k\in \mathbb N_0} \big\langle f, \psi_{\ell, j,k}\big\rangle \phi_{\ell, j,k}\eqno(2.8)$$

\parindent=0mm \vspace{.1in}
in the $L^2$-sense, then we say that ${\cal F}( \Phi)$ is a dual wavelet frame of ${\mathcal F}(\Psi)$ (and vice versa) or we simply say that (${\mathcal F}(\Psi), {\cal F}( \Phi)$) is a {\it pair of dual wavelet frames}.

\parindent=8mm \vspace{.2in}
In order to obtain a fast wavelet frame transform, tight wavelet frames are generally derived from refinable functions via multiresolution analysis.  We say that $\varphi\in L^2(K)$ is a {\it refinable function}, if it satisfies an equation of the type

~~~~~$$\varphi(x)=\sqrt q\sum_{k\in \mathbb N_0} c_k\varphi\big({\mathfrak p}^{-1}x- u(k)\big),\eqno(2.9)$$

\parindent=0mm \vspace{.1in}
where $c_k$ are complex coefficients. The functional equation $(2.9)$ is known as the {\it refinement equation}. Applying the Fourier transform, we can write this equation as

\parindent=0mm \vspace{.1in}
$$ \hat\varphi\left( \xi \right)=m_{0}({\mathfrak p}\xi)\hat\varphi({\mathfrak p}\xi),\eqno(2.10) $$

\parindent=0mm \vspace{.0in}
where $$m_{0}(\xi)=\dfrac{1}{\sqrt q}\sum_{k\in \mathbb N_0} c_{k}\overline{\chi_k(\xi) }.\eqno(2.11) $$

\parindent=0mm \vspace{.1in}
 Further, it is  proved  that a function $\varphi\in L^2(K)$  generates an MRA in $L^2(K)$ if and only if

$$ \sum_{k\in \mathbb N_0}\left| \hat\varphi\big(\xi- u(k)\big)\right|^2=1,~\text{for}~a.e.~\xi \in {\mathfrak D},\quad \lim_{j \rightarrow \infty}\left|\hat\varphi({\mathfrak p}^j\xi)\right|=1,~~\text{for}~~a.e.~\xi\in  K.\eqno(2.12)$$

\parindent=0mm \vspace{.2in}
Let the refinable function $\varphi\in L^2(K)$  generates an MRA $\left\{ V_j\right\}_{j\in\mathbb Z}$ of $L^2(K)$ and  $\Psi := \left\{\psi_1, \dots , \psi_L\right\}\subset  V_1$, then

$$ \psi_{\ell}\left(x\right)=\sqrt q\sum_{k \in \mathbb N_0}d^{\ell}_k\,\varphi\big( {\mathfrak p}^{-1}x- u(k)\big),~~\ell=1,\dots,L.\eqno(2.13) $$

Taking the Fourier transform for both sides of (2.13) gives

$$ \hat\psi_{\ell}\left( \xi \right)=m_{\ell}( {\mathfrak p}\xi)\hat\varphi( {\mathfrak p}\xi),$$

where
$$m_{\ell}(\xi)=\dfrac{1}{\sqrt q}\sum_{k \in \mathbb N_0} d_{k}^\ell\overline{\chi_k(\xi) },~~\ell=1,\dots,L \eqno(2.14) $$

\parindent=0mm \vspace{.1in}
are the {\it framelet symbols} or {\it wavelet masks}. With $m_\ell(\xi), \ell=0,1,\dots,L, L \ge q-1$, as the wavelet masks, we form the {\it modulation matrix} as

$${\mathcal H}(\xi)=\left(\begin{array}{cccc}
m_0(\xi)&m_0\big(\xi+ {\mathfrak p}u(1)\big)&\dots&m_0\big(\xi+ {\mathfrak p}u(q-1)\big)\\
m_1(\xi)&m_1\big(\xi+ {\mathfrak p}u(1)\big)&\dots&m_1\big(\xi+ {\mathfrak p}u(q-1)\big)\\
\vdots&\vdots&\ddots&\vdots\\
m_L(\xi)&m_L\big(\xi+ {\mathfrak p}u(1)\big)&\dots&m_L\big(\xi+ {\mathfrak p}u(q-1)\big)\\
\end{array}\right).\eqno(2.15)$$

\parindent=0mm \vspace{.2in}
The so-called unitary extension principle (UEP) provides a sufficient condition on $\Psi=\left\{\psi_1, \dots , \psi_L\right\}$ such that the  wavelet system ${\mathcal F}(\Psi)$ given by (2.5) constitutes a tight frame for $L^2(K)$. It is well known that in order to apply the UEP to derive wavelet tight frame from a given refinable function, the corresponding refinement mask must satisfy

$$\sum_{k=0}^{q-1}\left|m_{0}\big(\xi+{\mathfrak p}u(k)\big)\right|^{2}\le 1,\quad \xi\in K.\eqno(2.16)$$

\parindent=0mm \vspace{.1in}
Recently, Shah [14] has given a general procedure for the construction of   tight wavelet frames generated by the Walsh polynomials  using unitary extension principles as:

\parindent=0mm \vspace{.2in}
{\bf{Theorem 2.1.}} {\it Let  $\varphi(x)$ be a compactly supported refinable function and $\hat\varphi(0)=1$. Then, the wavelet system ${\mathcal F}(\Psi)$ given by (2.5) constitutes a Parseval frame in $L^2(K)$ provided the matrix ${\mathcal H}(\xi)$ as defined in (2.15) satisfies }

$${\mathcal H}(\xi){\mathcal H^*}(\xi)=I_q,\quad \text{for}~a.e.~\xi \in \sigma (V_0)\eqno(2.17)$$

\parindent=0mm \vspace{.1in}
{\it where $\sigma (V_0):=\big\{ \xi\in {\mathfrak D}:  \sum_{k\in \mathbb N_0}|\hat\varphi\big(\xi+u(k)\big)|^2\ne0  \big\}.$  }

\parindent=0mm \vspace{.2in}
{\bf{3.  Orthogonal Wavelet Frames on Local Fields }}

\parindent=0mm \vspace{.2in}
Motivated and inspired by the construction of tight wavelet frames generated by the Walsh polynomials [14] using the machinery of unitary extension principles.  In this section, we shall first derive the complete characterization of tight wavelet frames generated by the wavelet masks by means of their polyphase components.

\parindent=8mm \vspace{.2in}
The {\it polyphase representation} of the  refinement mask $m_0(\xi)$ can be derived as

$$\begin{array}{rcl}
m_0(\xi)&=&\displaystyle\dfrac{1}{\sqrt q}\sum_{k\in \mathbb N_0}h_k\; \overline{\chi_k(\xi)}\\\
&=&\displaystyle\dfrac{1}{\sqrt q}\sum_{r=0}^{q-1}\sum_{k\in \mathbb N_0}h_{u(r)+qk}\, \overline{\chi_{u(r)+qk}(\xi)}\\\
&=&\displaystyle\dfrac{1}{\sqrt q}\sum_{r=0}^{q-1}\overline{\chi_{u(r)}(\xi)}\sum_{k\in \mathbb N_0}h_{u(r)+qk}\, \overline{\chi_{qk}(\xi)}\\\
&=&\displaystyle\dfrac{1}{\sqrt q}\sum_{r=0}^{q-1}\overline{\chi_{u(r)}(\xi)}f_r^0\big(\overline{\chi(q\xi)}\big)
\end{array}$$

\parindent=0mm \vspace{.0in}
where

$$f_r^0(x)=\sum_{k\in \mathbb N_0}h_{u(r)+qk}\overline{x(k)},\quad r=0,1,\dots, q-1,\;x\in K.\eqno(3.1)$$

\parindent=0mm \vspace{.1in}
Similarly, the framelet symbols $m_\ell(\xi), \ell=1,2,\dots,L$, in defined in (2.14) can be splitted into polyphase components as

$$m_\ell(\xi)=\dfrac{1}{\sqrt q}\sum_{r=0}^{q-1}\overline{\chi_{u(r)}(\xi)}f_r^\ell\big(\overline{\chi(q\xi)}\big),\eqno(3.2)$$
where
$$f_r^\ell(\xi)=\sum_{k\in \mathbb N_0}h_{u(r)+qk}^\ell\overline{x(k)},\quad r=0,1,\dots, q-1,\;x\in K.\eqno(3.3)$$

\parindent=8mm \vspace{.1in}
With the polyphase components as defined in (3.1) and (3.3), we formulate the {\it polyphase matrix} $\Gamma\big(\overline{\chi(q\xi)}\big)$ as:

$$\Gamma\big(\overline{\chi(q\xi)}\big)=\left(\begin{array}{cccc}
f_0^0\big(\overline{\chi(q\xi)}\big)&f_0^1\big(\overline{\chi(q\xi)}\big)&\dots& f_0^L\big(\overline{\chi(q\xi)}\big)\\\
f_1^0\big(\overline{\chi(q\xi)}\big)&f_1^1\big(\overline{\chi(q\xi)}\big)&\dots &f_1^L\big(\overline{\chi(q\xi)}\big)\\
\vdots&\vdots&\ddots&\vdots\\
f_{q-1}^0\big(\overline{\chi(q\xi)}\big)&f_{q-1}^1\big(\overline{\chi(q\xi)}\big)&\dots& f_{q-1}^L\big(\overline{\chi(q\xi)}\big)
\end{array}\right).\eqno(3.4)$$

\parindent=0mm \vspace{.1in}
For convenience let $\overline{\chi(q\xi)}=\zeta$. The polyphase matrix  is called a {\it unitary matrix} if

 $${\Gamma}({\mathfrak p}^{-1}\zeta){\Gamma}^*({\mathfrak p}^{-1}\zeta)=I_q,\quad a.e.\,\xi\in {\mathfrak D}\eqno(3.5)$$

which is equivalent to
$$\sum_{\ell=0}^{L}\overline{f_r^\ell(\zeta)}f_{r^\prime}^\ell(\zeta)=\delta_{r,r^\prime} \Leftrightarrow \sum_{\ell=1}^{L}\overline{f_{r^\prime}^\ell(\zeta)}f_r^\ell(\zeta)=\delta_{r,r^\prime}-\overline{f_r^0(\zeta)}f_{r^\prime}^0(\zeta),\; 0\le r, r^\prime\le q-1. \eqno(3.6)$$

The following theorem shows that a unitary polyphase matrix leads to a tight wavelet frame on local fields of positive characteristic.

\parindent=0mm \vspace{.2in}
{\bf{Theorem 3.1.}} {\it  Suppose that the refinable function $\varphi$ and the framelet symbols $m_{0},m_{1},\dots, m_{L}$ satisfy equations (2.9)-(2.14). Furthermore, if the polyphase matrix $\Gamma(\zeta)$ given by (3.4) satisfy UEP condition (3.5), then the wavelet system ${\mathcal F}(\Psi)$ given by (2.5)  constitutes a tight frame for $L^2(K)$.}

\parindent=0mm \vspace{.2in}
{\bf{Proof.}} By Parseval's  formula, we have
\begin{align*}
\displaystyle\sum_{\ell=1}^{L}\sum_{ j \in \mathbb Z}\sum_{ k \in \mathbb N_0}\left|\big\langle f, \,\psi_{j,k}^\ell\big\rangle\right|^2&=
\sum_{\ell=1}^{L}\sum_{ j \in \mathbb Z}\sum_{ k \in \mathbb N_0}\left|\Big\langle f, \,q^{j/2}\psi^\ell\big({\mathfrak p}^{-j}x-u(k)\big)\Big\rangle\right|^2\\\
&=\sum_{\ell=1}^{L}\sum_{ j \in \mathbb Z}\displaystyle\sum_{ k \in \mathbb N_0}
\left| \left\langle \hat f, \,q^{j/2} \hat \psi^\ell \big({\mathfrak p}^{-j}\xi\big) \chi_{{\mathfrak p}^j}(\xi)\right\rangle\right|^2\\\
&=\sum_{\ell=1}^{L}\sum_{j \in \mathbb Z}q^j\sum_{ k \in \mathbb N_0}\left|\left\langle \hat f\big({\mathfrak p}^{-j}\xi \big)\overline{ \hat \psi^\ell(\xi)},\,\chi(\xi)\right\rangle\right|^2\\\
&=\sum_{\ell=1}^{L}\sum_{j \in \mathbb Z}q^j\int_K \left| \hat f({\mathfrak p}^{-j}\xi)\right|^2\left|\hat \psi^\ell(\xi)\right|^2 d\xi.\tag{3.7}
\end{align*}

\parindent=0mm \vspace{.1in}
Implementing the polyphase component formula (3.3) of wavelet masks  $m_\ell(\xi),\ell=1,\dots,L$, we can write
\begin{align*}
\sum_{\ell=1}^{L}\left|\hat \psi^\ell(\xi)\right|^2&=\sum_{\ell=1}^{L}\big| m_\ell({\mathfrak p}\xi) \hat \varphi({\mathfrak p}\xi)\big|^2\\\
&=\sum_{\ell=1}^{L}\overline{m_\ell({\mathfrak p}\xi)}\; \overline{\hat \varphi({\mathfrak p}\xi)}\,m_\ell({\mathfrak p}\xi)\hat \varphi({\mathfrak p}\xi)\\\
&= \overline{\hat \varphi({\mathfrak p}\xi)}\,\sum_{\ell=1}^{L}\overline{\left \{\dfrac{1}{\sqrt q}\sum_{r=0}^{q-1}\overline{\chi_{u(r)}({\mathfrak p}\xi)} f_r^\ell(\zeta)\right\}}\left \{\dfrac{1}{\sqrt q}\sum_{r^\prime=0}^{q-1}\overline{ \chi_{u({r^\prime})}({\mathfrak p}\xi)} f_{r^\prime}^\ell(\zeta)\right\} \hat \varphi({\mathfrak p}\xi)\\\
&= \overline{\hat \varphi({\mathfrak p}\xi)}\,\dfrac{1}{ q}\,\sum_{r=0}^{q-1}\sum_{r^\prime=0}^{q-1}\chi_{u(r)-u({r^\prime})}({\mathfrak p}\xi)\left \{\sum_{\ell=1}^{L} \overline{f_r^\ell(\zeta)}\, f_{r^\prime}^\ell(\zeta)\right\} \hat \varphi({\mathfrak p}\xi).
\end{align*}
Since the polyphase matrix ${\mathcal P}(\zeta)$  is unitary, which is equivalent to (3.6), the above expression reduces to
\begin{align*}
\sum_{\ell=1}^{L}\left|\hat \psi^\ell(\xi)\right|^2&=\overline{\hat \varphi({\mathfrak p}\xi)}\,\dfrac{1}{ q}\sum_{r=0}^{q-1}\sum_{r^\prime=0}^{q-1}\chi_{u(r)-u({r^\prime})}({\mathfrak p}\xi)\Big[\delta _{r, r^\prime}-f_r^0(\zeta)\, f_{r^\prime}^0(\zeta) \Big]\hat \varphi({\mathfrak p}\xi)\\\
\nonumber&=\overline{\hat \varphi({\mathfrak p}\xi)}\,\hat \varphi({\mathfrak p}\xi)-\overline{\hat \varphi({\mathfrak p}\xi)}\,\dfrac{1}{ q}\,\sum_{r=0}^{q-1}\sum_{r^\prime=0}^{q-1}\chi_{u(r)-u({r^\prime})}({\mathfrak p}\xi) \overline{f_r^0(\zeta)}\, f_{r^\prime}^0(\zeta)\hat \varphi({\mathfrak p}\xi)\\\
&=\left|\hat \varphi({\mathfrak p}\xi)\right|^2-\overline{\hat \varphi({\mathfrak p}\xi)}\,\overline{m_0({\mathfrak p}\xi)}\,m_0({\mathfrak p}\xi)\hat \varphi({\mathfrak p}\xi)\\\
&=\left|\hat \varphi({\mathfrak p}\xi)\right|^2-\left|m_0({\mathfrak p}\xi)\hat \varphi({\mathfrak p}\xi)\right|^2\\\
&=\left|\hat \varphi({\mathfrak p}\xi)\right|^2-\left|\hat \varphi(\xi)\right|^2.\tag{3.8}
\end{align*}

\parindent=0mm \vspace{.1in}
By substituting equation (3.8) in (3.7), we obtain
\begin{align*}
\sum_{\ell=1}^{L}\sum_{ j \in \mathbb Z}\sum_{ k \in \mathbb N_0}\left|\left\langle f,\, \psi_{j,k}^\ell\right\rangle\right|^2 &=\sum_{ j \in \mathbb Z} q^j \int_K \left|\hat f({\mathfrak p}^{-j}\xi)\right|^2\left\{\big|\hat \varphi({\mathfrak p}\xi)\big|^2-\big|\hat \varphi(\xi)\big|^2\right\} d \xi \\\
&=\int_K \left|\hat f(\xi)\right|^2\sum_{ j \in \mathbb Z}\left\{\left|\hat \varphi\big({\mathfrak p}^{j+1}\xi\big)\right|^2-\left|\hat \varphi\big({\mathfrak p}^j\xi\big)\right|^2\right\} d \xi.\tag{3.9}
\end{align*}

\parindent=0mm \vspace{.1in}
Using the assumption (2.12), the summand in the above expression can be written as

\begin{eqnarray*}
\sum_{ j \in \mathbb Z}\left\{\left|\hat \varphi\big({\mathfrak p}^{j+1}\xi\big)\right|^2-\left|\hat \varphi\big({\mathfrak p}^j\xi\big)\right|^2\right\} d \xi&=&\lim_{j \to \infty}\left|\hat \varphi\big({\mathfrak p}^{j+1}\xi\big)\right|^2-\displaystyle\lim_{j \to  -\infty}\left|\hat \varphi\big({\mathfrak p}^j\xi\big)\right|^2\\\
&=&\lim_{j \to \infty}\left|\hat \varphi\big({\mathfrak p}^j\xi\big)\right|^2-\lim_{j \to  \infty}\left|\hat \varphi\big({\mathfrak p}^{-j}\xi\big)\right|^2\\\
&=& \left|\hat \varphi(0)\right|^2-\lim_{j \to  \infty}\left|\hat \varphi\big({\mathfrak p}^{-j}\xi\big)\right|^2\\\
&=&1.
\end{eqnarray*}
By  using the above estimate in equation (3.9), we have
\begin{equation*}\sum_{ j \in \mathbb Z}\sum_{ k \in \mathbb N_0}\sum_{\ell=1}^{L}\left|\left\langle f, \,\psi_{j,k}^\ell\right\rangle\right|^2
=\int_K \left|\hat f(\xi)\right|^2d\xi=\left\|\hat f\right\|_2^2=\big\| f\big\|_2^2.\end{equation*}

\parindent=0mm \vspace{.1in}
This completes the proof of the theorem.\quad \fbox\\

\parindent=8mm \vspace{.2in}
The orthogonality of a pair of wavelet is guaranteed if the following extra conditions are imposed. Given a collection of wavelet masks ${\bf M}=[m_{0}, m_{1},\dots, m_{L}]$. For $k=k=1,2,\dots, q-1$, consider the following matrices

$${\mathcal M(\xi)}=\left(\begin{array}{ccc}
m_0(\xi)&m_0\big(\xi+ {\mathfrak p}u(k)\big)\\
m_1(\xi)&m_1\big(\xi+ {\mathfrak p}u(k)\big)\\
\vdots&\vdots\\
m_L(\xi)&m_L\big(\xi+ {\mathfrak p}u(k)\big)
\end{array}\right),\;\;{\mathcal M_0(\xi)}=\left(\begin{array}{ccc}
m_1(\xi)&m_1\big(\xi+ {\mathfrak p}u(k)\big)\\
m_2(\xi)&m_2\big(\xi+ {\mathfrak p}u(k)\big)\\
\vdots&\vdots\\
m_L(\xi)&m_L\big(\xi+ {\mathfrak p}u(k)\big)
\end{array}\right).\eqno(3.10)$$

Let there be another wavelet frame whose wavelet masks are given by $\tilde m_{0}, \tilde m_{1},\dots, \tilde m_{L}.$ Denoting the matrices as in (3.10) for these wavelet masks by $\tilde{\mathcal M}(\xi)$ and $\tilde{\mathcal M}_{0}(\xi)$, respectively. With the above definitions, we present an algorithm for the construction of arbitrarily many orthogonal wavelet tight frames generated by the wavelet masks.

\parindent=0mm \vspace{.2in}
{\bf{Theorem 3.2.}} {\it Suppose that $\varphi$ and $\tilde \varphi$ are the refinable functions that satisfy UEP. Let the corresponding filters be $m_\ell, \tilde m_\ell,\,\ell=0,1,\dots,L$. Let the matrices ${\mathcal M(\xi)},{\mathcal M_0(\xi)}, \tilde{\mathcal M}(\xi)$ and $\tilde{ \mathcal M_0}(\xi)$ be as defined in (3.10). For all $k=1,2,\dots,q-1$, suppose the following matrix equations hold}

$${\mathcal M^*(\xi)}{\mathcal M(\xi)}=I_2,\quad \tilde{\mathcal M}^*(\xi)\tilde{ \mathcal M}(\xi)=I_2,\quad\text{and}\quad {\mathcal M_0(\xi)}\tilde{ \mathcal M_0}(\xi)=0.\eqno(3.11)$$

\parindent=0mm \vspace{.1in}
{\it For $1\le \ell \le L$, let $\hat\psi_{\ell}\left( \xi \right)=m_{\ell}( {\mathfrak p}\xi)\hat\varphi( {\mathfrak p}\xi)$ and $\hat{\phi}_{\ell}\left( \xi \right)=m_{\ell}( {\mathfrak p}\xi)\hat{\tilde\varphi}( {\mathfrak p}\xi)$ be the corresponding dual. Then $\{\psi_1, \psi_2, \dots, \psi_L\}$ and $\{\phi_1, \phi_2, \dots, \phi_L\}$ generate orthogonal Parseval wavelet frames i.e., the systems $X(\Psi)$ and $X(\Phi)$ are orthogonal}.

\parindent=0mm \vspace{.2in}
{\it Proof.} From the Unitary Extension Principle, it follows that $\{\psi_1, \psi_2, \dots, \psi_L\}$ and $\{\phi_1, \phi_2,\\ \dots, \phi_L\}$ generate Parseval wavelet frames. It only remains to prove the orthogonality. For each $\ell$, by Holder's inequality and by virtue of the fact that $\psi_\ell$ and $\phi_\ell$ generate Bessel sequences, we have

$$\sum_{j\in \mathbb Z}\left|\hat\psi_\ell({\mathfrak p}^{-j}\xi)\overline{\hat{\phi}_\ell({\mathfrak p}^{-j}\xi)}\right|\le \sum_{j\in \mathbb Z}\left|\hat\psi_\ell({\mathfrak p}^{-j}\xi)\right|^2\sum_{j\in \mathbb Z}\left|\hat{\phi}_\ell({\mathfrak p}^{-j}\xi)\right|^2<\infty.\eqno(3.12)$$

\parindent=0mm \vspace{.1in}
Thus, the order of summation can be changed. With this, by equation (3.11), we have

$$\begin{array}{rcl}
\displaystyle\sum_{\ell=1}^L\sum_{j\in \mathbb Z}\hat\psi_\ell({\mathfrak p}^{-j}\xi)\overline{\hat{\phi}_\ell({\mathfrak p}^{-j}\xi)}&=&\displaystyle\sum_{\ell=1}^L\sum_{j\in \mathbb Z}m_\ell({\mathfrak p}^{1-j}\xi)\hat \varphi({\mathfrak p}^{1-j}\xi)\overline{\hat{\tilde m}_\ell({\mathfrak p}^{1-j}\xi)\hat{\tilde \varphi}({\mathfrak p}^{1-j}\xi)}\\
&=&\displaystyle\sum_{j\in \mathbb Z}\hat \varphi({\mathfrak p}^{1-j}\xi)\overline{\hat{\tilde \varphi}({\mathfrak p}^{1-j}\xi)}\sum_{\ell=1}^Lm_\ell({\mathfrak p}^{1-j}\xi)\overline{\hat{\tilde m}_\ell({\mathfrak p}^{1-j}\xi)}\\
&=&0,
\end{array}$$

\parindent=0mm \vspace{.1in}
holds for almost every $\xi \in \mathbb R^+$. Likewise, for $k\in \mathbb Z^+\setminus p\mathbb Z^+$, again by (3.11), we obtain

$$\begin{array}{lcr}
\displaystyle\sum_{\ell=1}^L\sum_{j=0}^\infty\hat\psi_\ell({\mathfrak p}^{-j}\xi)\overline{\hat{\phi}_\ell\big({\mathfrak p}^{-j}\big(\xi+ u(s)\big)\big)}&&\\
\qquad\qquad=\displaystyle\sum_{\ell=1}^L\sum_{j=0}^\infty m_\ell({\mathfrak p}^{1-j}\xi)\hat \varphi({\mathfrak p}^{1-j}\xi)\overline{\hat{\tilde m}_\ell\big({\mathfrak p}^{1-j}\big(\xi+ u(s)\big)\big)\hat{\tilde \varphi}\big({\mathfrak p}^{1-j}\big(\xi+ u(s)\big)\big)}&&\\
\qquad\qquad=\displaystyle\sum_{j=0}^\infty\hat \varphi({\mathfrak p}^{1-j}\xi)\overline{\hat{\tilde \varphi}\big({\mathfrak p}^{1-j}\big(\xi+ u(s)\big)\big)}\sum_{\ell=1}^Lm_\ell({\mathfrak p}^{1-j}\xi)\overline{\hat{\tilde m}_\ell\big({\mathfrak p}^{1-j}\big(\xi+ u(s)\big)\big)}&&\\
\qquad\qquad=0,
\end{array}$$

\parindent=0mm \vspace{.1in}
This completes the proof of the theorem.\quad \fbox\\

\parindent=8mm \vspace{.2in}
Next, we briefly describe how to obtain a pair of compactly supported orthogonal tight frames from a given compactly supported tight frame system ${\cal F}(\Psi)$ constructed via the UEP. More precisely, we construct a pair of orthogonal wavelet frames generated by the wavelet masks for the space $L^2(K)$ with  slightly different approach as described in Theorem 3.2.

\parindent=8mm \vspace{.2in}
Let $A$ be a $2L\times 2L$ paraunitary matrix. Partition $A=(A_1: A_2)$ where $A_1$ and $A_2$ are the first and last $L$ columns of $A$. Let $B$ and $C$ be the matrices

$$B=\left(\begin{array}{ccc}
1&&0\\
0&&A_1
\end{array}\right),\qquad C=\left(\begin{array}{ccc}
1&&0\\
0&&A_2
\end{array}\right).$$

\parindent=0mm \vspace{.1in}
With $B$ and $C$ in hand, we construct new polyphase matrices as $\Gamma_1=B\Gamma,\; \Gamma_2=C\Gamma$. The new polyphase matrix $\Gamma_1$ looks like

$$\begin{array}{rcl}
\Gamma_1&=&\left(\begin{array}{cccc}
1&0&\cdots&0\\
0&a_{1,1}&\cdots&a_{1,L}\\
\vdots&\vdots&\ddots&\vdots\\
0&a_{2L,1}&\cdots&a_{2L,L}
\end{array}\right)\left(\begin{array}{cccc}
f_0^0(\zeta)&f_1^0(\zeta)&\dots&f_{q-1}^0(\zeta)\\
f_0^1(\zeta)&f_1^1(\zeta)&\dots&f_{q-1}^1(\zeta)\\
\vdots&\vdots&\ddots&\vdots\\
f_0^L(\zeta)&f_1^L(\zeta)&\dots&f_{q-1}^L(\zeta)\\
\end{array}\right)\\\\
&=&\left(\begin{array}{cccc}
f_0^0(\zeta)&f_1^0(\zeta)&\dots& f_{q-1}^0(\zeta)\\\\
\displaystyle\sum_{\ell=1}^La_{1,L}f^\ell_0(\zeta)&\displaystyle\sum_{\ell=1}^La_{1,\ell}f^\ell_1(\zeta)&\dots &\displaystyle\sum_{\ell=1}^La_{1,\ell}f^\ell_{q-1}(\zeta)\\
\vdots&\vdots&\ddots&\vdots\\
\displaystyle\sum_{\ell=1}^La_{2L,\ell}f^\ell_0(\zeta)&\displaystyle\sum_{\ell=1}^La_{2L,\ell}f^\ell_1(\zeta)&\dots &\displaystyle\sum_{\ell=1}^La_{2L,\ell}f^\ell_{q-1}(\zeta)
\end{array}\right).
\end{array}$$

\parindent=0mm \vspace{.1in}
It is easy to verify that both the matrices $\Gamma_1$ and $\Gamma_2$  constructed above are unitary. Moreover, under this algorithm the scaling function does not change. Therefore, for $k=1,2,\dots,2L,$ the new wavelet masks $G_k(\xi)$ are given by

\begin{align*}
G_k(\xi)&=\dfrac{1}{\sqrt q}\displaystyle\sum_{r=0}^{q-1}\overline{\chi_{u(r)}(\xi)}\sum_{\ell=1}^La_{k,\ell}({\mathfrak p}^{-1}\xi)f_r^\ell(\zeta)\\
&=\sum_{\ell=1}^La_{k,\ell}({\mathfrak p}^{-1}\xi)\dfrac{1}{\sqrt q}\sum_{r=0}^{q-1}\overline{\chi_{u(r)}(\xi)}f_r^\ell(\zeta)\\
&=\sum_{\ell=1}^La_{k,\ell}({\mathfrak p}^{-1}\xi)m_\ell(\xi)\tag{3.13}
\end{align*}

\parindent=0mm \vspace{.1in}
Likewise one obtains $\tilde G_k(\xi)$ as

$$\tilde G_k(\xi)=\sum_{\ell=L+1}^{2L}a_{k,\ell}({\mathfrak p}^{-1}\xi)\tilde m_\ell(\xi).\eqno(3.14)$$

\parindent=0mm \vspace{.2in}
Let ${\mathcal M}(\xi)$ and $\tilde{\mathcal M}(\xi)$ be as in equation (3.10).  Then, ${\mathcal M^*(\xi)}{\mathcal M(\xi)}=I_2,\,\tilde{\mathcal M}^*(\xi)\tilde{ \mathcal M}(\xi)=I_2,$ as both the matrices ${\mathcal M}$ and $\tilde{\mathcal M}$ consist of the columns of the modulation matrices. This satisfies one of the conditions of Theorem 3.2.

\parindent=0mm \vspace{.2in}
{\bf{Lemma 3.3.}} {\it Let  $ {\mathcal M_0}(\xi)$ and $\tilde{ \mathcal M_0}(\xi)$ be the matrices of wavelet masks as in Theorem 3.2. Then }
$$ {\mathcal M_0}(\xi)^*  \tilde{ \mathcal M_0}(\xi)=0.\eqno(3.15)$$

\parindent=0mm \vspace{.2in}
{\it Proof.} Since the entries of the matrix $A$ are polynomials, so they are periodic in each components. Therefore, we have

$$A^\ell_r\big({\mathfrak p}^{-1}(\xi+{\mathfrak p}u(k))\big)=A^\ell_r\big({\mathfrak p}^{-1}\xi+u(k)\big)=A^\ell_r\big({\mathfrak p}^{-1}\xi\big).$$

\parindent=0mm \vspace{.1in}
Hence, (3.13) and (3.14) can be expressed as:

$$G_k\big(\xi+ {\mathfrak p}u(\ell)\big)=\sum_{\ell=1}^La_{k,\ell}({\mathfrak p}^{-1}\xi)m_\ell\big(\xi+ {\mathfrak p}u(\ell)\big), \quad \tilde G_k\big(\xi+ {\mathfrak p}u(\ell)\big)=\sum_{\ell=L+1}^{2L}a_{k,\ell}({\mathfrak p}^{-1}\xi)\tilde m_\ell\big(\xi+ {\mathfrak p}u(\ell)\big).$$

\parindent=0mm \vspace{.1in}
Thus, the matrix $ {\mathcal M_0}(\xi)$ in (2.17) becomes

$$\begin{array}{rcl}
{\mathcal M_0}(\xi)&=&\left(
\begin{array}{ccc}
\displaystyle\sum_{\ell=1}^La_{1,L}({\mathfrak p}^{-1}\xi)m_\ell(\xi)&&\displaystyle\sum_{\ell=1}^La_{1,\ell}({\mathfrak p}^{-1}\xi)m_\ell\big(\xi+ {\mathfrak p}u(\ell)\big)\\
\displaystyle\sum_{\ell=1}^La_{2,L}({\mathfrak p}^{-1}\xi)m_\ell(\xi)&&\displaystyle\sum_{\ell=1}^La_{2,\ell}({\mathfrak p}^{-1}\xi)m_\ell\big(\xi+ {\mathfrak p}u(\ell)\big)\\
\vdots&&\vdots\\
\displaystyle\sum_{\ell=1}^La_{L,\ell}({\mathfrak p}^{-1}\xi)m_\ell(\xi)&&\displaystyle\sum_{\ell=1}^La_{L,\ell}({\mathfrak p}^{-1}\xi)m_\ell\big(\xi+ {\mathfrak p}u(\ell)\big)
\end{array}\right)\\\\
&=&\left(\begin{array}{cccc}
a_{1,1}({\mathfrak p}^{-1}\xi)&a_{1,2}({\mathfrak p}^{-1}\xi)&\cdots&a_{1,L}({\mathfrak p}^{-1}\xi)\\
a_{2,1}({\mathfrak p}^{-1}\xi)&a_{2,1}({\mathfrak p}^{-1}\xi)&\cdots&a_{2,L}({\mathfrak p}^{-1}\xi)\\
\vdots&\vdots&\ddots&\vdots\\
a_{L,1}({\mathfrak p}^{-1}\xi)&a_{L,2}({\mathfrak p}^{-1}\xi)&\cdots&a_{L,L}({\mathfrak p}^{-1}\xi)
\end{array}\right)\left(\begin{array}{ccc}
m_1(\xi)&&m_1(\xi+ {\mathfrak p}u(\ell))\\
m_2(\xi)&&m_2(\xi+ {\mathfrak p}u(\ell))\\
\vdots&&\vdots\\
m_L(\xi)&&m_L(\xi+ {\mathfrak p}u(\ell))
\end{array}\right).
\end{array}$$

\parindent=0mm \vspace{.1in}
The corresponding dual matrix $\tilde{ \mathcal M_0}(\xi)$ is obtained similarly. Therefore, using the fact that the matrix $A$ is paraunitary, (3.15) holds.\qquad \fbox\\

\parindent=8mm \vspace{.2in}
For $k=1,2,\dots,2L$, define the wavelet system

$$\hat\psi^*_k\left({\mathfrak p}^{-1}\xi \right)=G_k( \xi)\hat\varphi(\xi), \quad \hat{\phi}^*_k\left({\mathfrak p}^{-1}\xi \right)=\tilde G_k( \xi)\hat{\tilde\varphi}(\xi)$$

\parindent=0mm \vspace{.1in}
and let $\Psi^*=\big\{\psi_1^*, \psi_2^*, \dots, \psi_{2L}^*\big\}$ and $\Phi^*=\big\{\phi_1^*, \phi_2^*, \dots, \phi_{2L}^*\big\}$.

\parindent=0mm \vspace{.2in}
{\bf{Theorem 3.4.}} {\it The wavelet systems $X(\Psi^*)$ and $X(\Phi^*)$ generated by $\{\psi_1^*, \psi_2^*, \dots, \psi_{2L}^*\}$ and $\{\phi_1^*, \phi_2^*, \dots, \phi_{2L}^*\}$ are a pair of orthogonal wavelet frames for $L^2(K)$.}

\parindent=0mm \vspace{.2in}
{\it Proof.} The proof of the theorem follows immediately from Theorem 3.1, Lemma 3.3 and the fact that the matrices $\Gamma_{1}$ and $\Gamma_{2}$ are  unitary.\qquad \fbox\\

\parindent=8mm \vspace{.2in}
The following result show the relationship between a pair of orthogonal MRA based wavelet frames.

\parindent=0mm \vspace{.2in}
{\bf{Theorem 3.5.}} {\it Suppose that ${\cal F}(\Psi)$ and ${\cal F}(\Phi)$ are a pair of orthogonal MRA wavelet frames for $L^2(K)$. If $P(\Psi)=P(\Phi)$ and there exists functions $h,g \in L^2(\mathbb R^+)$ such that $\Psi^g:=\{\psi_1^g, \psi_2^g, \dots, \psi_L^g\}$ and $\Phi^h:=\{\phi_1^h ,\phi_2^h, \dots, \phi_L^h\}$ are wavelet frames, where $\psi_\ell^g$ and $\phi_\ell^h$ are defined by $\hat \psi_\ell^g(\xi)=\hat \psi_\ell(\xi)\hat g(\xi),\,\hat \phi_\ell^h(\xi)=\hat \phi_\ell(\xi)\hat h(\xi), \,1\le \ell \le L$, respectively. Then, ${\cal F}(\Psi^g)$ and ${\cal F}(\Phi^h)$ are a pair of orthogonal wavelet frames for $L^2(K)$.}

\parindent=0mm \vspace{.2in}
{\it Proof.} Suppose that ${\cal F}(\Psi)$ and ${\cal F}(\Phi)$ are wavelet frames for $L^2(K)$ and $P(\Psi)=P(\Phi)$. Then, by the property of MRA based wavelet frames, for any $n\ne m \in \mathbb Z$, we have $P({\mathfrak p}^{-m}\Psi) \perp P({\mathfrak p}^{-n}\Phi)$. Therefore, for all $f_1\in P(\Psi)$, we have
\begin{align*}
 0&=Pf_1(x)\\\
&=\sum_{\ell=1}^L\sum_{j \in \mathbb Z}\sum_{k \in \mathbb N_0}\big\langle f_1(x), \psi_\ell\big({\mathfrak p}^{-j}x- u(k)\big)\big\rangle\phi_\ell\big({\mathfrak p}^{-j}x- u(k)\big)\\\
&=\sum_{\ell=1}^L\sum_{k \in \mathbb N_0}C_{ f_1, k}^\ell \phi_\ell\big(x- u(k)\big),\tag{3.16}
\end{align*}

\parindent=0mm \vspace{.1in}
where $C_{ f_1, k}^\ell =\big\langle f_1(x), \psi_\ell\big(x-u(k)\big)\big\rangle$. For any $f \in L^2(K)$, we define $f=f_1+f_2$, where $f_1\in P(\Psi),\, f_2\in L^2(K)\setminus P(\Psi)$, then, $\left\langle f_1, f_2\right\rangle=0.$ With this, we get

$$Pf_2(x)=\sum_{\ell=1}^L\sum_{k \in \mathbb N_0}\big\langle f_2(x), \psi_\ell\big(x-u(k)\big)\big\rangle\phi_\ell\big(x- u(k)\big)=\sum_{\ell=1}^L\sum_{k \in \mathbb N_0}C_{ f_2, k}^\ell\phi_\ell\big(x- u(k)\big)=0.\eqno(3.17)$$

\parindent=0mm \vspace{.1in}
By combining (3.16) and (3.17), we conclude that

$$Pf(x)=Pf_1(x)+Pf_2(x)=0.\eqno(3.18)$$

\parindent=0mm \vspace{.1in}
Since  $\hat \phi_\ell^g(\xi)=\hat \phi_\ell(\xi)\hat g(\xi)$ and $Pf(x)=0$, we have
\begin{align*}
0&=\widehat {Pf(x)}\\\
&=\sum_{\ell=1}^L\sum_{k \in \mathbb N_0}C_{ f, k}^\ell\hat \phi_\ell(\xi)\overline{\chi_k(\xi)}\\\
&=\sum_{\ell=1}^L\sum_{k \in \mathbb N_0}C_{ f_1, k}^\ell\hat\phi_\ell(\xi)\overline{\chi_k(\xi)}+\sum_{\ell=1}^L\sum_{k \in \mathbb N_0}C_{ f_2, k}^\ell\hat\phi_\ell(\xi)\overline{\chi_k(\xi)}\\\
&=\hat g(\xi)\left(\displaystyle\sum_{\ell=1}^L\sum_{k \in \mathbb N_0}C_{ f_1, k}^\ell\hat\phi_\ell(\xi)\overline{\chi_k(\xi)}+\sum_{\ell=1}^L\sum_{k \in \mathbb N_0}C_{ f_2, k}^\ell\hat\phi_\ell(\xi)\overline{\chi_k(\xi)}\right)\\\
&=\sum_{\ell=1}^L\sum_{k \in \mathbb N_0}C_{ f_1, k}^\ell\hat\phi_\ell^g(\xi)\overline{\chi_k(\xi)}+\sum_{\ell=1}^L\sum_{k \in \mathbb N_0}C_{ f_2, k}^\ell\hat\phi_\ell^g(\xi)\overline{\chi_k(\xi)}\\\
&=\sum_{\ell=1}^L\sum_{k \in \mathbb N_0}C_{ f, k}^\ell\hat \phi_\ell^g(\xi)\overline{\chi_k(\xi)}.\tag{3.19}
\end{align*}

\parindent=0mm \vspace{.1in}
Applying Inverse Fourier transform to the system (3.19), we get

$$ 0=Pf(x)=\sum_{\ell=1}^L\sum_{k \in \mathbb N_0}C_{ f, k}^\ell \phi_\ell^g(\xi)\overline{\chi_k(\xi)}.\eqno(3.20)$$

\parindent=0mm \vspace{.1in}
From the above equality, we deduce that
\begin{align*}\sum_{\ell=1}^L\sum_{k \in \mathbb N_0}C_{ f, k}^\ell\big\langle f(x), \phi_\ell^g\big(x- u(k)\big)\big\rangle&=\left\langle f(x),  \sum_{\ell=1}^L\sum_{k \in \mathbb N_0}C_{ f, k}^\ell\phi_\ell^g\big(x- u(k)\big)\right\rangle\\\
&=\big\langle f(x), 0\big\rangle\\\
&=\left\langle f(x), \sum_{\ell=1}^L\sum_{k \in \mathbb N_0}\big\langle f(x), \phi_\ell^g\big(x- u(k)\big)\big\rangle\psi_\ell\big(x- u(k)\big)\right\rangle\tag{3.21}
\end{align*}

\parindent=0mm \vspace{.1in}
Thus, we have

$$\sum_{\ell=1}^L\sum_{k \in \mathbb N_0}\big\langle f(x), \phi_\ell^g\big(x- u(k)\big)\big\rangle\psi_\ell\big(x- u(k)\big)=0.\eqno(3.22)$$

\parindent=0mm \vspace{.1in}
In a similar manner, we can show that

$$\sum_{\ell=1}^L\sum_{k \in \mathbb N_0}\big\langle f(x), \phi_\ell^g\big(x- u(k)\big)\big\rangle\psi_\ell^h\big(x- u(k)\big)=0.\eqno(3.23)$$

\parindent=0mm \vspace{.1in}
For any $j \in \mathbb Z$, we have

$$\begin{array}{lcr}
\displaystyle\sum_{\ell=1}^L\sum_{k \in \mathbb N_0}\big\langle f(x), \phi_\ell^g\big({\mathfrak p}^{-j}x- u(k)\big)\big\rangle\psi_\ell^h\big({\mathfrak p}^{-j}x- u(k)\big)&&\\
\qquad\qquad=q^j\displaystyle\sum_{\ell=1}^L\sum_{k \in \mathbb N_0}\left\langle f({\mathfrak p}^jx), \phi_\ell^g\big(x- u(k)\big)\right\rangle\psi_\ell^h\big(x- u(k)\big)
&&\\
\qquad\qquad=0.\qquad\qquad\qquad\qquad\qquad\qquad\qquad\qquad\qquad\qquad\qquad\qquad\qquad\qquad\qquad(3.24)
\end{array}$$

\parindent=0mm \vspace{.1in}
Putting everything together,  we conclude that

$$\sum_{\ell=1}^L\sum_{j \in \mathbb Z}\sum_{k \in \mathbb N_0}\big\langle f(x), \phi_\ell^g\big({\mathfrak p}^{-j}x- u(k)\big)\big\rangle\psi_\ell^h\big({\mathfrak p}^{-j}x- u(k)\big)=0.$$

\parindent=0mm \vspace{.1in}
Hence, ${\cal F}(\Psi^g)$ and ${\cal F}(\Phi^h)$ constitutes a pair of orthogonal wavelet frames generated by Walsh polynomials for $L^2(K)$.\quad \fbox\\

\parindent=8mm \vspace{.2in}
The following theorem describes a general construction algorithm for orthogonal wavelet tight frames for local fields of positive characteristic.

\parindent=0mm \vspace{.2in}
{\bf{Theorem 3.6.}} {\it Suppose $A(\xi)$ is an $L\times L$ paraunitary matrix with integral periodic entries $a_{\ell, r}(\xi)$ and let $A_r(\xi)$ denotes the $r$th column.  Let $m_{0},m_{1},\dots, m_{L}$ be the Walsh polynomials (masks) given by (2.11) and (2.14) such that ${\bf M}^*(\xi){\bf M}(\xi)=I_2$, where ${\bf M}=[m_0(\xi), m_1(\xi),\dots,m_L(\xi)]$ is the combined mask of the wavelet masks, and  let the wavelet system ${\mathcal F}(\Psi)$ forms a normalized wavelet frame for $L^2(K)$. For $r=1,2,\dots,L$, define  new wavelet masks via }

$$\begin{array}{rcl}
\left(\begin{array}{c}
\eta_{1,1}^r(\xi)\\
\eta_{1,2}^r(\xi)\\
\vdots\\
\eta_{1,L}^r(\xi)\\
\eta_{2,1}^r(\xi)\\
\vdots\\
\eta_{2,L}^r(\xi)\\
\vdots\\
\eta_{N,1}^r(\xi)\\
\eta_{N,2}^r(\xi)\\
\vdots\\
\eta_{N,L}^r(\xi)
\end{array}\right)&=&\left(\begin{array}{c}
A_r(\xi)m_1(\xi)\\\\
A_r(\xi)m_2(\xi)\\
\vdots\\
A_r(\xi)m_N(\xi)
\end{array}\right).
\end{array}\eqno(3.25)$$

\parindent=0mm \vspace{.1in}
{\it Then, for $r=1,2,\dots,L$, the affine systems generated by $\Psi^r=\{\psi_{n,\ell}^r:1\le n\le N,\,1\le \ell \le L\}$, where

$$\hat\psi_{n,\ell}^r({\mathfrak p}^{-1}\xi)=\eta_{n,\ell}^r(\xi)\hat\varphi(\xi),\eqno(3.26)$$

\parindent=0mm \vspace{.1in}
are tight wavelet frames and are pairwise orthogonal.}

\parindent=0mm \vspace{.2in}
{\it Proof.} We first prove that the systems ${\cal F}(\Psi^m), 1\le m\le L$ are tight wavelet frames for $L^2(K)$. To do so, we first consider

$${\mathcal M}_r=\left[m_0(\xi), \eta_{1,1}^r(\xi), \dots, \eta_{1,L}^r(\xi),\eta_{2,1}^r(\xi), \dots,\eta_{2,L}^r(\xi), \dots, \eta_{N,1}^r(\xi),  \dots,\eta_{N,L}^r(\xi)\right].$$

 Then we define ${\mathcal M}_r(\xi)$ according to (3.10) as

$$\begin{array}{rcl}
{\mathcal M}_r(\xi)&=&\left(\begin{array}{ccc}
m_0(\xi)&&m_0\big(\xi+ {\mathfrak p}u(k)\big)\\
\eta_{1,1}^r(\xi)&&\eta_{1,1}^r\big(\xi+ {\mathfrak p}u(k)\big)\\
\eta_{1,2}^r(\xi)&&\eta_{1,2}^r\big(\xi+ {\mathfrak p}u(k)\big)\\
\vdots&&\vdots\\
\eta_{1,L}^r(\xi)&&\eta_{1,L}^r\big(\xi+ {\mathfrak p}u(k)\big)\\
\eta_{2,1}^r(\xi)&&\eta_{2,1}^r\big(\xi+ {\mathfrak p}u(k)\big)\\
\vdots&&\vdots\\
\eta_{2,L}^r(\xi)&&\eta_{2,L}^r\big(\xi+ {\mathfrak p}u(k)\big)\\
\vdots&&\vdots\\
\eta_{N,1}^r(\xi)&&\eta_{N,1}^r\big(\xi+ {\mathfrak p}u(k)\big)\\
\eta_{N,2}^r(\xi)&&\eta_{N,2}^r\big(\xi+ {\mathfrak p}u(k)\big)\\
\vdots&&\vdots\\
\eta_{N,L}^r(\xi)&&\eta_{N,L}^r\big(\xi+ {\mathfrak p}u(k)\big)
\end{array}\right),
\end{array}\eqno(3.27)$$

\parindent=0mm \vspace{.1in}
for $k=1,2,\dots, q-1$. Then, ${\mathcal M}^*_r(\xi){\mathcal M}_r(\xi)$ is a $2 \times 2$  matrix. Next, we examine the entries of ${\mathcal M}^*_r(\xi){\mathcal M}_r(\xi)$ individually. Since the columns of $A(\xi)$ have length 1, it follows that

$$\begin{array}{rcl}
\left[{\mathcal M}^*_r(\xi){\mathcal M}_r(\xi)\right]_{1,1}&=&|m_0(\xi)|^2+\displaystyle\sum_{\ell=1}^L\sum_{n=1}^N\left|a_{\ell,r}(\xi)m_n(\xi)\right|^2\\
&=&|m_0(\xi)|^2+\displaystyle\sum_{\ell=1}^L|a_{\ell,r}(\xi)|^2\sum_{n=1}^N|m_n(\xi)|^2\\
&=&|m_0(\xi)|^2+\displaystyle\sum_{n=1}^N|m_n(\xi)|^2\\
&=&1.
\end{array}$$

\parindent=0mm \vspace{.1in}
Similarly,

$$\begin{array}{rcl}
\left[{\mathcal M}^*_r(\xi){\mathcal M}_r(\xi)\right]_{2,2}&=&\left|m_0\big(\xi+{\mathfrak p}u(k)\big)\right|^2+\displaystyle\sum_{\ell=1}^L\sum_{n=1}^N\left|a_{\ell,r}\big(\xi+{\mathfrak p}u(k)\big)m_n\big(\xi+{\mathfrak p}u(k)\big)\right|^2\\
&=&\left|m_0\big(\xi+{\mathfrak p}u(k)\big)\right|^2+\displaystyle\sum_{\ell=1}^L\sum_{n=1}^N\left|a_{\ell,r}\big(\xi+{\mathfrak p}u(k)\big)\right|^2|m_n\big(\xi+{\mathfrak p}u(k)\big)|^2\\
&=&\left|m_0\big(\xi+{\mathfrak p}u(k)\big)\right|^2+\displaystyle\sum_{n=1}^N|m_n\big(\xi+{\mathfrak p}u(k)\big)|^2\\
&=&1.
\end{array}$$

\parindent=0mm \vspace{.1in}
Using the fact that ${\mathcal M}^*(\xi){\mathcal M}(\xi)=I_2$ and that the entries of $A(\xi)$ are integral periodic, we have

$$\begin{array}{rcl}
\left[{\mathcal M}^*_r(\xi){\mathcal M}_r(\xi)\right]_{1,2}&=&m_0\big(\xi+{\mathfrak p}u(k)\big)\overline{m_0(\xi)}+\displaystyle\sum_{\ell=1}^L\sum_{n=1}^Na_{\ell,r}\big(\xi+{\mathfrak p}u(k)\big)m_n\big(\xi+{\mathfrak p}u(k)\big)
\overline{a_{\ell,r}(\xi)m_n(\xi)}\\
&=&m_0\big(\xi+{\mathfrak p}u(k)\big)\overline{m_0(\xi)}+\displaystyle\sum_{\ell=1}^L\sum_{n=1}^N|a_{\ell,r}(\xi)|^2\overline{m_n(\xi)}m_n\big(\xi+{\mathfrak p}u(k)\big)\\
&=&m_0\big(\xi+{\mathfrak p}u(k)\big)\overline{m_0(\xi)}+\displaystyle\sum_{n=1}^N\overline{m_n(\xi)}m_n\big(\xi+{\mathfrak p}u(k)\big)\\
&=&0.
\end{array}$$

\parindent=0mm \vspace{.1in}
By the conjugate symmetry of ${\mathcal M}^*_r(\xi){\mathcal M}_r(\xi)$, the entry (2,1) must be zero. Thus

$${\mathcal M}^*_r(\xi){\mathcal M}_r(\xi)=I_2,\quad 1\le r \le L.\eqno(3.28)$$

\parindent=0mm \vspace{.1in}
Putting everything together, from Theorem 3.2, the wavelet systems ${\cal F}(\Psi^m)$ defined via (3.26) are tight wavelet frames for $L^2(K)$. It only remains to prove the orthogonality. According to equation (3.10),  for $1\le r \le L,$ we have

$${\mathcal M}^0_r(\xi)=\begin{array}{rcl}
\left(\begin{array}{ccc}
\eta_{1,1}^r(\xi)&\eta_{1,1}^r\big(\xi+ {\mathfrak p}u(k)\big)\\
\eta_{1,2}^r(\xi)&\eta_{1,2}^r\big(\xi+ {\mathfrak p}u(k)\big)\\
\vdots&\vdots\\
\eta_{1,L}^r(\xi)&\eta_{1,L}^r\big(\xi+ {\mathfrak p}u(k)\big)\\
\eta_{2,1}^r(\xi)&\eta_{2,1}^r\big(\xi+ {\mathfrak p}u(k)\big)\\
\vdots&\vdots\\
\eta_{2,L}^r(\xi)&\eta_{2,L}^r\big(\xi+ {\mathfrak p}u(k)\big)\\
\vdots&\vdots\\
\eta_{N,1}^r(\xi)&\eta_{N,1}^r\big(\xi+ {\mathfrak p}u(k)\big)\\
\eta_{N,2}^r(\xi)&\eta_{N,2}^r\big(\xi+ {\mathfrak p}u(k)\big)\\
\vdots&\vdots\\
\eta_{N,L}^r(\xi)&\eta_{N,L}^r\big(\xi+ {\mathfrak p}u(k)\big)
\end{array}\right)&=\left(\begin{array}{ccc}
A_r(\xi)m_1(\xi)&A_r(\xi + {\mathfrak p}u(k))m_1\big(\xi+ {\mathfrak p}u(k)\big)\\\\
A_r(\xi)m_2(\xi)&A_r(\xi + {\mathfrak p}u(k))m_2\big(\xi+ {\mathfrak p}u(k)\big)\\
\vdots&\vdots\\
A_r(\xi)m_N(\xi)&A_r(\xi + {\mathfrak p}u(k))m_N\big(\xi+ {\mathfrak p}u(k)\big)
\end{array}\right)
\end{array}\eqno(3.29)$$

\parindent=0mm \vspace{.1in}
If $1\le r\ne r^\prime\le L$, then

$$\begin{array}{lcr}
{\mathcal M}^0_r(\xi)^*{\mathcal M}^0_r(\xi)&&\\\\
\quad =\left(\begin{array}{ccc}
A_r(\xi)m_1(\xi)&&A_r(\eta)m_1(\eta)\\
A_r(\xi)m_2(\xi)&&A_r\big(\eta)m_2(\eta)\\
\vdots&&\vdots\\
A_r(\xi)m_N(\xi)&&A_r(\eta)m_N(\eta)
\end{array}\right)^*\left(\begin{array}{ccc}
A_{r^\prime}(\xi)m_1(\xi)&&A_{r^\prime}(\eta)m_1(\eta)\\
A_{r^\prime}(\xi)m_2(\xi)&&A_{r^\prime}(\eta)m_2(\eta)\\
\vdots&&\vdots\\
A_{r^\prime}(\xi)m_N(\xi)&&A_{r^\prime}(\eta)m_N(\eta)
\end{array}\right)&&\\\\
\quad=\left(\begin{array}{ccc}
A_r^*(\xi)A_n(\xi)\displaystyle\sum_{n=1}^N|m_n(\xi)|^2&&A_r^*(\xi)A_{r^\prime}(\eta)\displaystyle\sum_{n=1}^N\overline{m_n(\xi)}m_n(\eta)\\
A_r^*(\eta)A_{r^\prime}(\xi)\displaystyle\sum_{n=1}^N\overline{m_n(\eta)}m_n(\xi)&&A_r^*(\eta)A_{r^\prime}(\eta)\displaystyle\sum_{n=1}^N
\overline{m_n(\eta)}m_n(\eta)\\
\end{array}\right)&&\\\\
\quad=\left(\begin{array}{ccc}
A_r^*(\xi)A_{r^\prime}(\xi)\displaystyle\sum_{n=1}^N|m_n(\xi)|^2&&A_r^*(\xi)A_{r^\prime}(\xi)\displaystyle\sum_{n=1}^N\overline{m_n(\xi)}m_n(\eta)\\
A_r^*(\eta)A_{r^\prime}(\xi)\displaystyle\sum_{n=1}^N\overline{m_n(\eta)}m_n(\xi)&&A_r^*(\eta)A_{r^\prime}(\xi)\displaystyle\sum_{n=1}^N
\overline{m_n(\eta)}m_n(\eta)\\
\end{array}\right)&&\\
\quad=0,
\end{array}$$

\parindent=0mm \vspace{.1in}
where $\eta=\xi+\mathfrak p u(k)$. Here, we have used the fact that the product of the two matrices $A_m^*(\xi)A_{m^\prime}(\xi)=0$ by the
orthogonality of the columns of $A(\xi)$. Using Theorem 3.2, we get the desired result.\qquad \fbox\\

\parindent=0mm \vspace{.2in}
{\bf{References}}

\end{document}